\documentclass[11pt]{article}
\usepackage{amsfonts,epsf,amsmath,amssymb}
\usepackage{tikz}
\usepackage{cite}

\newtheorem{thm}{\bf Theorem}

\newtheorem{lemma}[thm]{\bf Lemma}

\newtheorem{cor}[thm]{\bf Corollary}
\newtheorem{proposition}[thm]{\bf Proposition}

\newenvironment{unnumbered}[1]{\trivlist \item [\hskip \labelsep {\bf
#1}]\ignorespaces\it}{\endtrivlist}

\newcommand{\proof}{\noindent{\bf Proof.\ }}
\newcommand{\qed}{\hfill $\square$ \medskip}

\newcommand{\gtd}{\gamma_{{\rm tg}}}

\newcommand{\D}{{Dominator }}
\newcommand{\St}{{Staller }}

\newcommand{\barU}{\overline{U}}

\let\oldenumerate\enumerate
\renewcommand{\enumerate}{
  \oldenumerate
  \setlength{\itemsep}{0pt}
  \setlength{\parskip}{0pt}
  \setlength{\parsep}{0pt}
}

\newcommand{\2}{ \vspace{0.2cm} }

\begin{document}

\title{\bf Game Total Domination Critical Graphs}

\author{
Michael A. Henning $^{a}$
\and
Sandi Klav\v zar $^{b,c,d}$
\and
Douglas F. Rall $^{e}$
}

\date{}

\maketitle

\begin{center}
$^a$ Department of Pure and Applied Mathematics \\ University of Johannesburg,
South Africa \\
{\tt mahenning@uj.ac.za}
\medskip

$^b$ Faculty of Mathematics and Physics, University of Ljubljana, Slovenia\\
$^c$ Faculty of Natural Sciences and Mathematics, University of Maribor, Slovenia\\
$^d$ Institute of Mathematics, Physics and Mechanics, Ljubljana, Slovenia\\
{\tt sandi.klavzar@fmf.uni-lj.si}
\medskip

$^e$ Department of Mathematics, Furman University, Greenville, SC, USA\\
{\tt doug.rall@furman.edu}
\end{center}

\begin{abstract}
In the total domination game played on a graph $G$, players Dominator and Staller alternately select vertices of $G$, as long as possible, such that each vertex chosen increases the number of vertices totally dominated. Dominator (Staller) wishes to minimize (maximize) the number of vertices selected. The game total domination number, $\gamma_{\rm tg}(G)$, of $G$ is the number of vertices chosen when Dominator starts the game and both players play optimally. If a vertex $v$ of $G$ is declared to be already totally dominated, then we denote this graph by $G|v$. In this paper the total domination game critical graphs are introduced as the graphs $G$ for which $\gamma_{\rm tg}(G|v) < \gamma_{\rm tg}(G)$ holds for every vertex $v$ in $G$. If $\gamma_{\rm tg}(G) = k$, then $G$ is called $k$-$\gamma_{\rm tg}$-critical. It is proved that the cycle $C_n$ is $\gamma_{{\rm tg}}$-critical if and only if $n\pmod 6 \in \{0,1,3\}$ and that  the path $P_n$ is $\gamma_{{\rm tg}}$-critical if and only if $n\pmod 6\in \{2,4\}$. $2$-$\gamma_{\rm tg}$-critical and $3$-$\gamma_{\rm tg}$-critical graphs are also characterized as well as $3$-$\gamma_{\rm tg}$-critical joins of graphs.
\end{abstract}

\noindent{\bf Keywords:} total domination game; game total domination number; critical graphs; paths and cycles

\medskip
\noindent{\bf AMS Subj. Class.:} 05C57, 05C69, 05C38

\newpage
\section{Introduction}

Just as total domination~\cite{heya-2013} is considered as the natural counterpart of the standard domination, the total domination game is studied in parallel to the domination game. The latter game was introduced in~\cite{brklra-2010} and extensively studied afterwards, papers~\cite{BDK-2017, BDK-2016, brklko-2013, Bu-2015, bujtas-2015, doko-2013, heki-2016, hl-2017, bill-2013, gk-2012, nadjafi-2016, sc-2016} form a selection of the developments on the game.

A vertex $u$ in a graph $G$ \emph{totally dominates} a vertex~$v$ if $u$ is adjacent to $v$ in $G$. A \emph{total dominating set} of $G$ is a set $S$ of vertices of $G$ such that every vertex of $G$  is totally dominated by a vertex in $S$.  The \emph{total domination game} was introduced in~\cite{hkr-2015} as follows. Given a graph $G$, two players, called \emph{Dominator} and \emph{Staller}, take turns choosing a vertex from $G$. Each vertex chosen must totally dominate at least one vertex not totally dominated by the set of vertices previously chosen. Such a chosen vertex is called a \emph{legal move}. The game ends when $G$ there is no legal move available. Dominator wishes to minimize the number of vertices selected, while the goal of Staller is just the opposite. If Dominator (Staller, resp.) has the first move, then we speak about a {\em D-game} ({\em S-game}, resp.). The \emph{game total domination number}, $\gtd(G)$, of $G$ is the number of moves played in the D-game when both players play optimally. The corresponding invariant for the S-game is denoted $\gtd'(G)$.

Adopting the terminology from~\cite{bill-2013}, a \emph{partially total dominated graph} is a graph together with a declaration that some vertices are already totally dominated; that is, they need not be totally dominated in the rest of the game. Given a graph $G$ and a subset $S$ of vertices of $G$, we denote by $G|S$ the partially total dominated graph in which the vertices of $S$ in $G$ are already totally dominated. If $S = \{v\}$ for some vertex $v$ in $G$, we simply write $G|v$. We use $\gtd(G|S)$ to denote the number of turns remaining in the D-game on $G|S$.

As already mentioned, the study of the total domination game was initiated in~\cite{hkr-2015}, where it was demonstrated that the two versions differ significantly. In~\cite{dh-2016} $\gtd$ as well as $\gtd'$ is determined for paths and cycles. Then, in~\cite{hkr-2017} it was proved that if $G$ is a graph of order $n$ in which every component contains at least three vertices, then $\gtd(G) \le 4n/5$. It was further conjectured that $3n/4$ is the correct upper bound. Significant progress on this conjecture was made in~\cite{BuHeTu16, hr-2016}. In~\cite{hr-2017} trees with equal total domination and game total domination number are characterized, while in~\cite{brhe-2017} it is proved that the game total domination problem is log-complete in PSPACE.

Critical graphs with respect to the domination game were introduced and studied in~\cite{buklko-2015}. In this paper we introduce a parallel concept for the total domination game as follows. A graph $G$ is \emph{total domination game critical}, abbreviated $\gtd$-critical, if $\gtd(G) > \gtd(G|v)$ holds for every $v \in V(G)$.  If $G$ is $\gtd$-critical and $\gtd(G) = k$, we say that $G$ is $k$-$\gtd$-critical. Upgrading the main results from~\cite{dh-2016} we characterize in Sections~\ref{sec:cycles} and~\ref{sec:paths} $\gtd$-critical cycles and paths, respectively. In the final section we  characterize $2$-$\gamma_{\rm tg}$-critical and $3$-$\gamma_{\rm tg}$-critical graphs as well as joins of graphs that are $3$-$\gamma_{\rm tg}$-critical.


In the rest of this section we recall some definitions and notation. The \emph{open neighborhood} of a vertex $v$ in $G$ is $N_G(v) = \{u \in V(G) \, | \, uv \in E(G)\}$ and its \emph{closed neighborhood } is $N_G[v] = \{v\} \cup N_G(v)$. If the graph $G$ is clear from the context, we simply write $N(v)$ and $N[v]$ rather than $N_G(v)$ and $N_G[v]$, respectively.  For a set $S \subseteq V(G)$, the subgraph induced by $S$ is denoted by $G[S]$. A vertex $u$ in a graph $G$ \emph{totally dominates} a vertex~$v$ if $v\in N_G(u)$. A \emph{total dominating set} of $G$ is a set $S$ of vertices of $G$ such that every vertex of $G$  is totally dominated by a vertex in $S$. A \emph{dominating vertex} of a graph $G$ is a vertex adjacent to every other vertex of $G$. For notation and graph theory terminology not defined herein, we in general follow~\cite{heya-2013}. Finally, we use the standard notation $[k] = \{1,\ldots,k\}$.

\section{Cycles}
\label{sec:cycles}

In this section we prove the following result.

\begin{thm}
\label{t:critcycle}
For $n \ge 3$, the cycle $C_n$ is $\gtd$-critical if and only if $n \pmod 6 \in \{0,1,3\}$.
\end{thm}

We first recall that the game total domination number of a cycle is determined in~\cite{dh-2016}.

\begin{thm}{\rm (\cite{dh-2016})}
 \label{t:cycle_thm}
For $n \ge 3$,
\[
\gtd(C_n) =
 \left\{
 \begin{array}{cl}
 \lfloor\frac{2n+1}{3} \rfloor -1   & \mbox{when $n \equiv 4 \pmod 6$} \2 \\
 \lfloor\frac{2n+1}{3} \rfloor  & \mbox{otherwise.} \end{array}
\right.
\]
\end{thm}

Following the notation from~\cite{dh-2016}, we define a \emph{run} in a partially total dominated cycle to be a maximal sequence of consecutive totally dominated vertices that contain at least two vertices, and an \emph{anti-run} as a maximal sequence of consecutive vertices none of which are totally dominated and that contain at least two vertices.

Suppose first that $n \pmod 6 \in \{0,1,3\}$ and consider the cycle $G \cong C_n$ given by $v_1 v_2 \ldots v_nv_1$. We show that $G$ is $\gtd$-critical. Let $v$ be an arbitrary vertex of $G$ and consider the D-game played on $G|v$. If $n = 3$, then $\gtd(G|v) = 1$ since \D plays as his first move the vertex~$v$. Hence, we may assume that $n \ge 6$. By symmetry we may assume that $v = v_1$. At each point in the game we denote by $U$ the set of vertices that are not playable, that is vertices already played or vertices not played but all of whose neighbors are already totally dominated. Each move of the game adds at least one vertex to the set $U$, namely the vertex played. \D plays as his first move the vertex $v_4$, thereby ensuring that the two vertices $v_2$ and $v_4$ are added to $U$ after Dominator's first move.

After Dominator's first move, his strategy is to guarantee that after each of Staller's moves, his answer together with her move add at least three vertices to $U$. His strategy is precisely that presented in~\cite{dh-2016}. For completeness, we describe this strategy. Consider a move of \St on some vertex $v_i$ where $i \in [k] \setminus \{4\}$. Thus, $v_{i-1}$ or $v_{i+1}$ is a new vertex totally dominated by the move on $v_i$. For notational convenience, we may assume that $v_{i+1}$ is a new vertex totally dominated by the move on $v_i$. Since $v_{i+1}$ was not totally dominated before Staller's move, the vertex $v_{i+2}$ was not played yet. Suppose $v_{i+4}$ was not played either. Then, $v_{i+3}$ is not yet totally dominated. In this case, \D can play $v_{i+4}$, and ensure that the three vertices $v_i$, $v_{i+2}$ and $v_{i+4}$ are added to $U$ after these two moves. On the other hand, suppose that $v_{i+4}$ was played before \St plays vertex $v_i$. Then already with Staller's move on $v_i$, both vertices $v_i$ and $v_{i+2}$ were added to $U$. So \D can play any legal move, adding at least one vertex to $U$, and thus ensuring that at least three vertices get added to $U$ after these two moves.

We next compute the bounds given by Dominator's strategy. As observed earlier, Dominator's opening move adds two vertices to $U$. Let $m$ be the number of moves played when the game finishes. Consider first the case when $m$ is odd, that is, when \D finishes the game. Then, Dominator's strategy ensures that $n = |U| \ge 2+\frac{3}{2}(m-1)$, implying that $m \le \frac{1}{3}(2n - 1)$. Now suppose that $m$ is even, and thus that \St made the final move. Let $v_i$ be a new vertex totally dominated by Staller's final move. Before \St played her final move, at least two moves were legal, namely $v_{i-1}$ and $v_{i+1}$. So with Staller's last move, she adds at least two vertices to $U$. Thus, in this case, $n = |U| \ge 2 + \frac{3}{2}(m-2) + 2$, and so $m \le \frac{2}{3}(n-1) < \frac{1}{3}(2n - 1)$. In both cases, $m \le \frac{1}{3}(2n - 1)$. Thus since $m$ is an integer, $m \le \lfloor \frac{1}{3}(2n - 1) \rfloor$. Therefore, noting that $n\pmod 6 \in\{0,1,3\}$, we get that $\gtd(C_n) \le \lfloor  \frac{1}{3}(2n - 1) \rfloor = \lfloor \frac{2n+1}{3} \rfloor - 1$. This is true for every vertex $v$ of $G$. By Theorem~\ref{t:cycle_thm}, $\gtd(G) = \lfloor\frac{2n+1}{3} \rfloor$. Thus, the graph $G$ is $\gtd$-critical. This proves the sufficiency.

Next we consider the necessity, which we prove using a contrapositive argument. Suppose that $n \pmod 6 \in \{2,4,5\}$. We show that the graph $G$ is not $\gtd$-critical. If $n = 4$ it is straightforward to check that $G$ is not $\gtd$-critical. Hence, we may assume that $n \ge 5$. Let $v$ be an arbitrary vertex of $G$. We describe a strategy for Staller that guarantees that she totally dominates exactly one new vertex on all of her moves, and then deduce that $\gtd(G|v) \ge \gtd(G)$. Staller's strategy follows her strategy presented in~\cite{dh-2016}. However for completeness, we describe her strategy in detail in order to enable us to compute the precise bounds.

Suppose it is Staller's turn to play, and denote by $A$ the subset of vertices already totally dominated. In particular, we note that $v \in A$. We may assume that $A$ does not contain all vertices or the game would be finished. If the cycle contains a run, i.e., a maximal sequence $v_i v_{i+1} \ldots v_{i+r}$ ($r\ge 1$) of consecutive vertices in $A$, then \St plays the run extremity $v_{i+r}$ and totally dominates only one new vertex, namely $v_{i+r+1}$. If the cycle contains an anti-run, i.e., a maximal sequence $v_jv_{j+1} \ldots v_{j+s}$  ($s\ge 1$) of consecutive vertices not in $A$, then \St plays on the anti-run extremity $v_{j+s}$ thereby totally dominating only one new vertex, namely $v_{j+s-1}$. Thus, Staller's strategy of choosing the extremity of a run or an anti-run whenever one exists ensures that only one new vertex is totally dominated whenever she plays such a move.

As observed in~\cite{dh-2016}, a situation when there is no run and no anti-run in the cycle may occur only if the cycle is of even length, in which case the cycle is bipartite with the set $A$ one of the partite sets. During the game, every move of \St therefore totally dominates exactly one new vertex, except possibly once when there is no run and no anti-run for \St to play on. In particular, if $n \equiv 5 \pmod 6$, then \St totally dominates exactly one new vertex on each of her moves. We note that every move of \D totally dominates at most two new vertices. Further, before the game starts the vertex $v$ is prescribed as totally dominated.

Suppose that $n \pmod 6\in \{2,4\}$ and that there is no run and no anti-run for \St to play on. In this case, the set $A$ is one of the partite sets and it is Staller's move. We note that every previous move of \St totally dominates one vertex. Recall that $v \in A$ and that \D made the first move. If $n \equiv 2 \pmod 6$, then  $|A| \equiv 1 \pmod 3$, while if $n \equiv 4 \pmod 6$, then  $|A| \equiv 2 \pmod 3$. Suppose that every move of \D totally dominated exactly two new vertices. If $|A| \equiv 1 \pmod 3$, it is necessarily Dominator's turn to move, a contradiction. If $|A| \equiv 2 \pmod 3$, then \D would necessarily have totally dominated exactly one new vertex on at least one of his moves, a contradiction. Therefore, \D totally dominated only one new vertex on at least one of his moves.

From our above observations, in this case when $n \pmod 6\in \{2,4,5\}$, either \St totally dominates exactly one new vertex on each of her moves or \St totally dominates exactly one new vertex on all except one of her moves and \D totally dominated one new vertex on at least one of his moves. Let $m$ be the number of moves played when the game finishes.

Assume first that $m$ is even. Then, \St made the final move and both players played $\frac{m}{2}$ moves. In this case, our earlier observations imply that \St totally dominates at most~$\frac{m}{2}$ vertices and \D totally dominates at most~$m$ vertices or \St totally dominates at most~$\frac{m}{2} + 1$ vertices and \D totally dominates at most~$m - 1$ vertices. In both cases, \St and \D together totally dominate at most $\frac{3m}{2}$ vertices. Recall that the vertex $v$ is already totally dominated before the game begins. Thus, $n \le \frac{3m}{2}+1$ and therefore $m \ge \lceil \frac{2(n-1)}{3}\rceil$. Observe that $\lceil \frac{2(n-1)}{3}\rceil = \lfloor\frac{2n+1}{3} \rfloor -1$  when $n \equiv 4 \pmod 6$ and $\lceil \frac{2(n-1)}{3}\rceil = \lfloor\frac{2n+1}{3} \rfloor$  when $n \pmod 6 \in \{2,5\}$. Thus, by Theorem~\ref{t:cycle_thm}, $m = \gtd(G|v) \ge \gtd(G)$.

Assume next that $m$ is odd. Then, \D made the final move and therefore played $\frac{m+1}{2}$ moves and \St played $\frac{m-1}{2}$ moves. In this case, our earlier observations imply that \St totally dominates at most~$\frac{m-1}{2}$ vertices and \D totally dominates at most~$m+1$ vertices or \St totally dominates at most~$\frac{m-1}{2} + 1$ vertices and \D totally dominates at most~$m$ vertices. In both cases, \St and \D together totally dominate at most~$\frac{3m+1}{2}$ vertices. Recall that the vertex $v$ is already totally dominated before the game begins. Thus, $n \le \frac{3(m+1)}{2}$, and therefore $m \ge \lceil \frac{2n-3}{3} \rceil$. If $n \pmod 6 \in \{2,5\}$, then $\lceil \frac{2n-3}{3} \rceil = \lfloor\frac{2n+1}{3} \rfloor$, while if $n \equiv 4 \pmod 6$, then $\lceil \frac{2n-3}{3} \rceil = \lfloor\frac{2n+1}{3} \rfloor - 1$. Thus in all cases, by Theorem~\ref{t:cycle_thm}, $m = \gtd(G|v) \ge \gtd(G)$. Hence, we have shown that if $n \pmod 6\in \{2,4,5\}$, then the graph $G$ is not $\gtd$-critical. This completes the proof of Theorem~\ref{t:critcycle}.

\section{Paths}
\label{sec:paths}

The main result of this section reads as follows.

\begin{thm}
\label{t:critpath}
For $n \ge 2$, the path $P_n$ is $\gtd$-critical if and only if $n \pmod 6\in \{2,4\}$.
\end{thm}

The game total domination number of a path is determined in~\cite{dh-2016}.

\begin{thm}{\rm (\cite{dh-2016})}
 \label{t:path_thm}
For $n \ge 2$,
\[
\gtd(P_n) =
 \left\{
 \begin{array}{cl}
 \lfloor \frac{2n}{3} \rfloor   & \mbox{when $n \equiv 5 \pmod 6$} \2 \\
 \lceil \frac{2n}{3}  \rceil  & \mbox{otherwise.} \end{array}
\right.
\]
\end{thm}

To prove Theorem~\ref{t:critpath}, we adopt the notation introduced in~\cite{dh-2016}. Suppose that we are playing on a partially total dominated path, $P_n$, given by $x_1y_1x_2y_2 \ldots x_qy_q$ if $n = 2q$ is even and $x_1y_1x_2y_2 \ldots x_qy_qx_{q+1}$ if $n = 2q+1$ is odd. As defined in~\cite{dh-2016}, let $X$ and $Y$ denote the two partite sets of $P_n$, where $X = \{x_1,x_2,\ldots,x_{\lceil n/2 \rceil}\}$ and $Y = \{y_1,y_2,\ldots,y_{\lfloor n/2 \rfloor}\}$, and partition the sets $X$ and $Y$ into subsets
\[
X = (X_1,X_2,\ldots,X_{\lceil \frac{n}{6} \rceil}) \hspace*{0.5cm} \mbox{and} \hspace*{0.5cm}
Y = (Y_1,Y_2,\ldots,Y_{\lceil \frac{n-1}{6} \rceil})
\]
where $X_i = \{x_{3i-2},x_{3i-1},x_{3i}\}$ for $i \in \left[\lceil \frac{n}{6} \rceil - 1\right]$, $Y_i = \{y_{3i-2},y_{3i-1},y_{3i}\}$ for $i \in \left[\lceil \frac{n-1}{6} \rceil - 1\right]$, and where
\[
X_{\lceil \frac{n}{6} \rceil} = X \setminus \left( \bigcup_{i=1}^{\lceil \frac{n}{6} \rceil - 1} X_i \right) \hspace*{0.5cm} \mbox{and} \hspace*{0.5cm}
Y_{\lceil \frac{n-1}{6} \rceil} = Y \setminus \left( \bigcup_{i=1}^{\lceil \frac{n-1}{6} \rceil - 1} Y_i \right).
\]

We note that $|X_{\lceil \frac{n}{6} \rceil}| \in [3]$ and $|Y_{\lceil \frac{n-1}{6} \rceil}| \in [3]$, while every other subset in the partition of $X$ or $Y$ has cardinality~$3$. Each subset of the partition with exactly one vertex (respectively, exactly two vertices) not currently totally dominated we call a \emph{single} (respectively, \emph{double}). Further, each subset of the partition with all three vertices not currently totally dominated we call a \emph{triple}. Let $t_1$, $t_2$ and $t_3$ denote the number of singles, doubles and triples, respectively.

For $m \ge 0$, at any particular stage in the game immediately after Dominator's $m$th move has been played on a path $P_n$, let
\[
f_n(m) = 2m - 1 + 2t_3 + \left\lceil \frac{3}{2}t_2 \right\rceil + t_1.
\]

An identical proof used to prove~\cite[Lemma~13]{dh-2016} yields the following result.

\begin{lemma}
\label{l:nonincreasing}
For $m \ge 0$, Dominator has a strategy to ensure that $f_n(m+1) \le f_n(m)$.
\end{lemma}

A similar proof used to prove~\cite[Lemma~16]{dh-2016} yields the following result. However for completeness, we provide a short proof of this result.

\begin{lemma}
\label{l:bound_Dom}
If $G \cong P_n$ and $v$ is an arbitrary vertex of $G$, then $\gtd(G|v) \le f_n(1)$.
\end{lemma}

\proof Suppose, firstly, that $\gtd(G|v)$ is odd, say $\gtd(G|v) = 2r + 1$ for some $r \ge 0$. The game is therefore completed immediately after Dominator plays his $(r+1)$st move, implying that after this move is played, $t_1 = t_2 = t_3 = 0$ and $f_n(r+1) = 2(r+1)-1 = \gtd(G|v)$. Thus, by Lemma~\ref{l:nonincreasing}, $\gtd(G|v) = f_n(r+1) \le f_n(1)$.
Suppose, next, that $\gtd(G|v)$ is even, say $\gtd(G|v) = 2r$. Thus, after Dominator plays his $r$th move, $t_1 = 1$ and $t_2 = t_3 = 0$ since one additional move (played by Staller) finishes the game. Thus, $f_n(r) = 2r - 1 + t_1 = 2r = \gtd(G|v)$. Thus, by Lemma~\ref{l:nonincreasing}, $\gtd(G|v) = f_n(r) \le f_n(1)$.
\qed

We are now in a position to prove the following lemmas.

\begin{lemma}
 \label{l:2mod6}
For $n \ge 2$ and $n \equiv 2 \pmod 6$, if $G \cong P_n$ and $v$ is an arbitrary vertex of $G$, then $\gtd(G|v) \le \gtd(G) - 1$.
\end{lemma}

\proof For $n \ge 2$ and $n \equiv 2 \pmod 6$, let $G \cong P_n$ and let $v$ be an arbitrary vertex of $G$. If $n = 2$, then $\gtd(G|v) = 1$ and $\gtd(G) = 2$. Hence let $n = 6k+2$ for some $k \ge 1$. We note that $X_{k+1} = \{x_{3k+1}\}$ and $X_i = \{x_{3i-2},x_{3i-1},x_{3i}\}$ for $i \in [k]$. Further, $Y_{k+1} = \{y_{3k+1}\}$ and $Y_i = \{y_{3i-2},y_{3i-1},y_{3i}\}$ for $i \in [k]$. By Theorem~\ref{t:path_thm}, $\gtd(G) = \lceil \frac{2n}{3}  \rceil = 4k+2$. Hence it suffices for us to show that $\gtd(G|v) \le 4k+1$.

If $v = y_{3k+1}$, then Dominator plays as his first move $y_{3k-1}$. If $v = x_{3k+1}$, then Dominator plays as his first move $y_{3k-1}$. If $v = y_{3i}$ for some $i \in [k]$, then Dominator plays as his first move $x_{3i-1}$. If $v = x_{3i}$ for some $i \in [k]$, then Dominator plays as his first move $y_{3i-2}$. If $v = y_{3k-1}$, then Dominator plays as his first move $x_{3k+1}$.  If $v = x_{3k-1}$, then Dominator plays as his first move $y_{3k}$. If $v = y_{3i-2}$ for some $i \in [k]$, then Dominator plays as his first move $x_{3i}$. If $v = x_{3i-2}$ for some $i \in [k]$, then Dominator plays as his first move $y_{3i-1}$. In all cases after Dominator's move, $t_1 = 2$, $t_2 = 0$ and $t_3 = 2k-1$, implying that $f_n(1) = 1 + 2(2k-1) + 2 = 4k+1$.

If $v \in \{x_{3i-1},y_{3i-1}\}$ for some $i \in [k-1]$ and $k \ge 2$, then Dominator plays as his first move $x_{3k+1}$. After Dominator's move, $t_1 = 1$, $t_2 = 2$ and $t_3 = 2k-2$, implying that $f_n(1) = 1 + 2(2k-2) + 3 + 1 = 4k+1$.

The above cases exhaust all possible choices for the vertex $v$. In all cases, $f_n(1) = 4k+1$. Thus, by Lemma~\ref{l:bound_Dom}, $\gtd(G|v) \le f_n(1) = 4k+1$.
\qed

\begin{lemma}
 \label{l:4mod6}
For $n \ge 4$ and $n \equiv 4 \pmod 6$, if $G \cong P_n$ and $v$ is an arbitrary vertex of $G$, then $\gtd(G|v) \le \gtd(G) - 1$.
\end{lemma}

\proof For $n \ge 4$ and $n \equiv 4 \pmod 6$, let $G \cong P_n$ and let $v$ be an arbitrary vertex of $G$. If $n = 4$, then $\gtd(G|v) = 2$ and $\gtd(G) = 3$. Hence assume that $n = 6k+4$ for some $k \ge 1$. We note that $X_{k+1} = \{x_{3k+1},x_{3k+2}\}$ and $X_i = \{x_{3i-2},x_{3i-1},x_{3i}\}$ for $i \in [k]$. Further, $Y_{k+1} = \{y_{3k+1},y_{3k+2}\}$ and $Y_i = \{y_{3i-2},y_{3i-1},y_{3i}\}$ for $i \in [k]$. By Theorem~\ref{t:path_thm}, $\gtd(G) = \lceil \frac{2n}{3}  \rceil = 4k+3$. Hence it suffices for us to show that $\gtd(G|v) \le 4k+2$.

If $v \in Y_{k+1}$, then Dominator plays as his first move $y_{3k+1}$. If $v \in X_{k+1}$, then Dominator plays as his first move $x_{3k+2}$. In all cases after Dominator's move, $t_1 = 1$, $t_2 = 0$ and $t_3 = 2k$, implying that $f_n(1) = 1 + 2(2k) + 1 = 4k+2$.

If $v = y_{3i}$ for some $i \in [k]$, then Dominator plays as his first move $x_{3i-1}$. If $v = x_{3i}$ for some $i \in [k]$, then Dominator plays as his first move $y_{3i-2}$. If $v \in \{ x_{3i-1}, y_{3i-1}\}$ for some $i \in [k]$, then Dominator plays as his first move $x_{3k+2}$. If $v = y_{3i-2}$ for some $i \in [k]$, then Dominator plays as his first move $x_{3i}$. If $v = x_{3i-2}$ for some $i \in [k]$, then Dominator plays as his first move $y_{3i-1}$. In all cases after Dominator's move, $t_1 = 0$, $t_2 = 2$ and $t_3 = 2k-1$, implying that $f_n(1) = 1 + 2(2k-1) + 3 = 4k+2$.

The above cases exhaust all possible choices for the vertex $v$. In all cases, $f_n(1) = 4k+2$. Thus, by Lemma~\ref{l:bound_Dom}, $\gtd(G|v) \le f_n(1) = 4k+2$.
\qed

As an immediate consequence of Lemmas~\ref{l:2mod6} and~\ref{l:4mod6}, we have the following result.

\begin{cor}
\label{c:2or4mod6}
For $n \ge 2$, if $n \pmod 6 \in \{2,4\}$, then the path $P_n$ is $\gtd$-critical.
\end{cor}

We prove next the following result.

\begin{lemma}
\label{l:not2or4mod6}
For $n \ge 3$, if $n \pmod 6 \not\in \{2,4\}$, then the path $P_n$ is not $\gtd$-critical.
\end{lemma}
\proof Let $G$ be the path $P_n$ given by $v_1v_2 \ldots v_n$ where $n \ge 3$ and $n \pmod 6 \in \{0,1,3,5\}$. It suffices for us to show that there exists a vertex $v$ in $G$ such that $\gtd(G|v) \ge \gtd(G)$. If $n = 3$, then letting $v = v_1$, we note that $\gtd(G|v) = 2 = \gtd(G)$. If $n = 5$, then letting $v = v_3$, we note that $\gtd(G|v) = 3 = \gtd(G)$. If $n = 6$, then letting $v = v_3$, we note that $\gtd(G|v) = 4 = \gtd(G)$. Hence in what follows we may assume that $n \ge 7$, for otherwise the desired result follows. Let $P$ and $Q$ be the partite sets of $G$ containing $v_1$ and $v_2$, respectively; that is, $P$ is exactly the set of all vertices $v_i$ with odd subscripts and $Q$ is exactly the set of all vertices $v_i$ with even subscripts.

\begin{unnumbered}{Claim~\ref{l:not2or4mod6}.1}
If $n \pmod 6 \in \{0, 1\}$, then $G$ is not $\gtd$-critical.
\end{unnumbered}
\proof Suppose that $n \pmod 6 \in \{0, 1\}$. If $n \equiv 0 \pmod 6$, then $n = 6k$ for some $k \ge 2$, and, by Theorem~\ref{t:path_thm}, $\gtd(G) = \lceil 2n/3 \rceil = 4k$. If $n \equiv 1 \pmod 6$, then $n = 6k + 1$ for some $k \ge 1$, and by Theorem~\ref{t:path_thm} $\gtd(G) = \lceil 2n/3 \rceil = 4k + 1$. We now select the vertex $v = v_4$ and show that $\gtd(G|v) \ge \gtd(G)$. Let $P_1 = \{v_1\}$ and $P_2 = P \setminus \{v_1,v_3,v_5\}$, and let $Q_1 = \{v_2\}$ and let $Q_2 = Q \setminus \{v_2,v_4\}$.  If $n = 6k$, we note that $|Q| = |P| = 3k$, $|Q_2| = 3k - 2$ and $|P_2| = 3(k-1)$. If $n = 6k+1$, we note that $|Q| = 3k$, $|Q_2| = 3k - 2$, $|P| = 3k+1$ and $|P_2| = 3k-2$. Staller's strategy is as follows.

Whenever \D plays a vertex in $Q$ (that totally dominates a new vertex in $P$) and there is a vertex in $P$ not yet totally dominated, then \St responds by playing a (legal) vertex in $Q$ at the extremity of a run or anti-run thereby totally dominating exactly one new vertex of $P$. We note that such a move of \St is always possible since in this case she plays a vertex in $Q$ in response to Dominator's move that plays a vertex in $Q$.

Whenever \D plays a vertex in $P$ that totally dominates a vertex in $Q_2$ and there is a vertex in $Q_2$ not yet totally dominated after Dominator's move, then \St responds by playing a (legal) vertex in $P$ at the extremity of a run or anti-run that totally dominates exactly one new vertex of $Q_2$. We note that such a move of \St is always possible since the vertex $v_4 \in Q$ is totally dominated.

Whenever \D plays the vertex $v_1$ or $v_3$ (which serves only to totally dominate the vertex $v_2 \in Q_1$), then \St responds as follows.
\begin{itemize}
\item[(i)] If there is a vertex in $Q_2$ not yet totally dominated, then she plays a vertex in $P$ at the extremity of a run or anti-run that totally dominates exactly one new vertex of $Q_2$.
\item[(ii)] If all vertices in $Q$ are totally dominated and at least one vertex in $P$ (but not all, for otherwise the game is finished) is totally dominated, then she plays a (legal) vertex in $Q$ at the extremity of a run or anti-run thereby totally dominating exactly one new vertex of $P$.
\item[(iii)] If all vertices in $Q$ are totally dominated and no vertex in $P$ is totally dominated, then she plays the vertex $v_4$.
\item[(iv)] Thereafter, \St follows her original strategy with the additional strategy that whenever \D plays a vertex in $Q_2$ and there is a vertex in $P_2$ that is not yet totally dominated, she plays a legal vertex in $Q_2$ at the extremity of a run or anti-run (such a move of \St totally dominates exactly one new vertex of $P_2$).
\end{itemize}

Let $m$ be the number of moves in the game played in $G|v$. Staller's strategy guarantees that she totally dominates exactly one new vertex on all of her moves, except for possible one move.

\medskip
\emph{Case~1. \St totally dominates two new vertices on one of her moves.} We note that this case can only occur when \D plays the vertex $v_1$ or $v_3$, thereby totally dominating exactly one new vertex (namely the vertex $v_2 \in Q_1$). Further, at this stage of the game all vertices in $Q$ are totally dominated and no vertex in $P$ is totally dominated. Thus, Staller's move that totally dominates two new vertices is the vertex $v_4$ in response to Dominator's move that plays $v_1$ or $v_3$.  Since $|Q_2| = 3(k-1) + 1$ and each of Staller's moves totally dominates at most one vertex of $Q_2$, at least~$2(k-1) + 1 = 2k-1$ moves were played to totally dominate all vertices of~$Q_2$. Therefore, at least~$2k$ moves are played to totally dominate all vertices of~$Q$.

After \St has played the vertex $v_4$, each of her subsequent moves totally dominates at most one vertex of $P_2$, hence at least~$\lceil 2|P_2|/3 \rceil$ moves are played to totally dominate all vertices of~$P_2$. Thus if $n = 6k$, then $|P_2| = 3(k-1)$ and at least $2(k-1)$ moves are played to totally dominate all vertices of~$P_2$, while if $n = 6k+1$, then $|P_2| = 3(k-1) + 1$ and at least $2k-1$ moves are played to totally dominate all vertices of~$P_2$. Further, one additional move is played to totally dominate the vertex $v_1$. Therefore, if $n = 6k$, then at least~$1 + 2(k-1) + 1 = 2k$ moves are played to totally dominate all vertices of~$P$, while if $n = 6k+1$, then at least~$1 + (2k-1) + 1 = 2k+1$ moves are played to totally dominate all vertices of~$P$.

Hence, if $n = 6k$, then $\gtd(G|v) = m \ge 2k + 2k = 4k = \gtd(G)$, while if $n = 6k + 1$, then $\gtd(G|v) = m \ge 2k + (2k+1) = 4k+1 = \gtd(G)$. In both cases, $\gtd(G|v) \ge \gtd(G)$.

\medskip
\emph{Case~2. \St totally dominates exactly one new vertex on all of her moves.} Staller's strategy forces \D to totally dominate the vertex $v_2 \in Q$ (by playing either $v_1$ or $v_3$). With this move, \D totally dominates exactly one new vertex in $Q$. Since $|Q_2| = 3(k-1) + 1$ and each of Staller's moves totally dominate at most one vertex of $Q_2$, at least~$2k-1$ moves are played to totally dominate all vertices of~$Q_2$. Therefore, at least~$2k$ moves are played to totally dominate all vertices of~$Q$. Recall that in this case, \St totally dominates exactly one new vertex on all of her moves. Thus if $n = 6k$, then $|P| = 3k$ and at least~$2k$ moves are played to totally dominate all vertices of~$P$, while if $n = 6k+1$, then $|P| = 3k+1$ and at least~$2k+1$ moves are played to totally dominate all vertices of~$P$.
Hence, if $n = 6k$, then $m \ge 2k + 2k = 4k$, while if $n = 6k + 1$, then $m \ge 2k + (2k+1) = 4k+1$. Thus, as in the proof of Case~1, we once again have that $\gtd(G|v) \ge \gtd(G)$. This completes the proof of Claim~\ref{l:not2or4mod6}.1.
\qed

\begin{unnumbered}{Claim~\ref{l:not2or4mod6}.2}
If $n \equiv 3 \pmod 6$, then $G$ is not $\gtd$-critical.
\end{unnumbered}
\proof Suppose that $n \equiv 3 \pmod 6$. Thus, $n = 6k+3$ for some $k \ge 1$. By Theorem~\ref{t:path_thm}, $\gtd(G) = \lceil 2n/3 \rceil = 4k+2$. We now select the vertex $v = v_1$ and show that $\gtd(G|v) \ge 4k+2$. Staller adopts the following simple strategy: she plays a vertex on the extremity of a run or an anti-run. We note that \St can always play according to her strategy since $n$ is odd and $v = v_1$. Thus, each of Staller's moves totally dominate exactly one new vertex. Let $m$ be the number of moves in the game played in $G|v$. Since in this case $|Q| = 3k+1$, at least $2k+1$ moves are played to totally dominate all vertices of~$Q$. Further, since $|P \setminus \{v_1\}| = 3k+1$, at least $2k+1$ moves are played to totally dominate all vertices of~$P$. Thus, $\gtd(G|v) = m \ge 4k+2 = \gtd(G)$.
\qed

\begin{unnumbered}{Claim~\ref{l:not2or4mod6}.3}
If $n \equiv 5 \pmod 6$, then $G$ is not $\gtd$-critical.
\end{unnumbered}
\proof Suppose that $n \equiv 5 \pmod 6$. Thus, $n = 6k+5$ for some $k \ge 1$. By Theorem~\ref{t:path_thm} $\gtd(G) = \lfloor 2n/3 \rfloor = 4k+3$. We choose $v = v_1$ and show that $\gtd(G|v) \ge 4k+3$. As in the proof of Claim~\ref{l:not2or4mod6}.2, \St plays a vertex on the extremity of a run or an anti-run. We note that \St can always play according to her strategy since $n$ is odd and $v = v_1$. Thus, each of Staller's moves totally dominate exactly one new vertex. Let $m$ be the number of moves in the game played in $G|v$.

Suppose firstly that $m$ is even. Then, \St made the final move, and so \D played $m/2$ moves and \St played $m/2$ moves. In this case, \D totally dominates at most~$m$ vertices, while \St totally dominates~$m/2$ vertices since as observed earlier she totally dominates exactly one new vertex on all of her moves. This implies that  $n \le |\{v\}| + m + m/2 = 1 + 3m/2$, or equivalently, $\gtd(G|v) = m \ge \lceil 2(n-1)/3 \rceil = 4k+3 = \gtd(G)$.

Suppose secondly that $m$ is odd. Then, \D made the final move, and so \D played $(m+1)/2$ moves and \St played $(m-1)/2$ moves. In this case, \D totally dominates at most~$m+1$ vertices, while \St totally dominates~$m-1$ vertices. This implies that  $n \le |\{v\}| + (m+1) + (m-1)/2 = 3(m+1)/2$, or equivalently, $\gtd(G|v) = m \ge \lceil (2n-3)/3 \rceil = 4k+3 = \gtd(G)$. Thus in both cases, $\gtd(G|v) = m \ge 4k+2 = \gtd(G)$.
\qed

The proof of Lemma~\ref{l:not2or4mod6} follows from Claim~\ref{l:not2or4mod6}.1,~\ref{l:not2or4mod6}.2 and~\ref{l:not2or4mod6}.3.
\qed

Theorem~\ref{t:critpath} immediately follows from Corollary~\ref{c:2or4mod6} and Lemma~\ref{l:not2or4mod6}.

\section{Characterization of $2$- and $3$-$\gtd$-critical graphs}
\label{sec:2-3-critical}

In this section we characterize  $2$-$\gtd$-critical graphs, $3$-$\gtd$-critical graphs, and joins of graphs that are $3$-$\gtd$-critical. The characterization of $2$-$\gtd$-critical graphs is simple.

\begin{proposition}
\label{prp:2-critical}
A graph is $2$-$\gtd$-critical if and only if it is a complete graph.
\end{proposition}

\proof Suppose that $G$ is a $2$-$\gtd$-critical graph of order~$n$. Thus, $\gtd(G) = 2$ and $\gtd(G|v) = 1$ for every vertex $v$ in $G$. In particular, $n \ge 2$ and $G$ has no isolated vertex. Let $v$ be an arbitrary vertex in $G$. Since $\gtd(G|v) = 1$, there is a vertex $v^*$ in $G$ that totally dominates every vertex of $G$ except for the vertex $v$. This is only possible if $v^* = v$ and $v$ is adjacent to every vertex in $G$ other than $v$ itself. Since $v$ is an arbitrary vertex of $G$, this implies that $G \cong K_n$. Conversely, if $G \cong K_n$ for some $n \ge 2$, then $\gtd(G) = 2$ and $\gtd(G|v) = 1$ for every vertex $v$ in $G$, implying that $G$ is  $2$-$\gtd$-critical.
\qed

To characterize $3$-$\gtd$-critical graphs some preparation is needed. Two vertices $u$ and $v$ of a graph $G$ are \emph{open twins} in $G$ if their open neighborhoods are the same; that is, $N(u) = N(v)$. A graph is \emph{open twin-free} if it contains no open twins. The proof of the following lemma is analogous to that presented in~\cite{buklko-2015} for domination game critical graphs. For completeness, we give a proof of this result.

\begin{lemma}
\label{l:lem1}
If $u$ and $v$ are open twins in $G$, then $\gtd(G) = \gtd(G|u) = \gtd(G|v)$.
\end{lemma}
\proof Suppose that the same game is played on $G$ and on $G|v$, that is, the same vertices are selected in both games. Then we claim that a move is legal in the game played on $G$ if and only if the move is legal in the game played on $G|v$. A legal move played in $G|v$ is clearly legal when played in $G$. If during the course of the game, a legal move $m^*$ played in $G$ is not legal in $G|v$, then $v$ would be the only newly totally dominated vertex in the game played on $G$, implying that the vertex $u$ was already totally dominated. However, any played vertex that totally dominates $u$ would also totally dominate $v$, since $u$ and $v$ are open twins. This contradicts the fact that immediately before the move $m^*$ is played in $G$ the vertex $v$ is not yet totally dominated. Therefore, at every stage of the game, a legal move in $G$ is legal in the game played on $G|v$. This proves the claim, from which the lemma follows immediately.
\qed

As an immediate consequence of Lemma~\ref{l:lem1}, we have the following result.

\begin{cor}
\label{c:co1}
Every $\gtd$-critical graph is open twin-free.
\end{cor}

In order to characterize the class of $3$-$\gtd$-critical graphs, we first establish some fundamental properties of such graphs, the first one being true of all $\gtd$-critical graphs.

\begin{lemma}
\label{l:lem2B}
If $G$ is a $\gtd$-critical graph and $v$ is any vertex of $G$, then no neighbor of $v$ is an optimal first move of \D in the D-game on $G|v$.
\end{lemma}

\proof Assume that $G$ is a $k$-$\gtd$-critical graph.  Let $v$ be an arbitrary vertex of $G$ and let $d_1$ be an optimal first move of \D in the D-game on $G|v$.  This
means that $k-1=\gtd(G|v)=1+\gtd'(G|(\{v\}\cup N(d_1))$. If $d_1 \in N(v)$, then \D could play $d_1$ as a first move in the game on $G$, which  implies that
\[\gtd(G)\le 1+\gtd'(G|N(d_1))=1+\gtd'(G|(\{v\}\cup N(d_1))=1+k-2=k-1\,,\]
a contradiction.
\qed

\begin{lemma}
\label{l:lem2}
If $G$ is a $3$-$\gtd$-critical graph, then $G$ has no dominating vertex.  Furthermore, for every vertex $w$ of $G$, $\gtd(G|w)=2$.
\end{lemma}

\proof Suppose that $G$ is a $3$-$\gtd$-critical graph. Thus, $\gtd(G) = 3$ and $\gtd(G|v) \le 2$ for every vertex $v$ in $G$.  If $G$ has a dominating vertex $v$, then
\D plays the vertex $v$ as his first move in the D-game played in $G$. The only legal move of \St in response to Dominator's move is to play a neighbor of $v$, implying that after her first move the game is complete. Hence, $\gtd(G) = 2$, a contradiction. Now, let $w$ be an arbitrary vertex of $G$.  Since $w$ is not a dominating vertex, it is clear that
$\gtd(G|w)\ge 2$.
\qed

We note that the first move $d_1$ of Dominator in the proof of Lemma~\ref{l:lem2B} could be $d_1=v$ as the example of the path $P_4$ with $v$ being a vertex of degree $2$ shows.

By Theorem~\ref{t:cycle_thm} and Theorem~\ref{t:critcycle}, no cycle is $3$-$\gtd$-critical. By Theorem~\ref{t:path_thm} and Theorem~\ref{t:critpath}, the path $P_4$ is the only $3$-$\gtd$-critical path. By Lemma~\ref{l:lem2}, a $3$-$\gtd$-critical graph has no dominating vertex. Thus, the only $3$-$\gtd$-critical graph of order~$4$ or less is the path $P_4$. More generally, we are now ready for the following characterization of $3$-$\gtd$-critical graphs. \newpage

\begin{thm}
\label{t:char3gtd}
Let $G$ be a graph of order~$n$ with no isolated vertex. The graph $G$ is $3$-$\gtd$-critical if and only if the following all hold. \\[-23.5pt]
\begin{enumerate}
\item The graph $G$ is open twin-free.
\item There is no dominating vertex in $G$.
\item For every vertex $v$ of $G$ of degree at most $n-3$, there exists a vertex $u$ of degree~$n-2$ that is not adjacent to $v$.
\end{enumerate}
\end{thm}

\proof Suppose firstly that properties (a), (b) and (c) in the statement of the theorem hold. In particular, this implies that $G$ has maximum degree~$\Delta(G) = n-2$. We show that $\gtd(G) = 3$. Let $d_1$ be an optimal first move of Dominator in the game played in $G$. Since $G$ has no dominating vertex, there is a vertex $x$ that is not adjacent to the vertex $d_1$ in $G$. Since $G$ is open twin-free, there is a vertex $y$ that is a neighbor of $d_1$ but not of $x$ or is a neighbor of $x$ but not of $d_1$. In both cases, \St can force at least three moves to be made to finish the game by playing the vertex $y$ on her first move. Thus, $\gtd(G) \ge 3$. On the other hand, if \D selects as his first move a vertex of degree~$n-2$, then after his move exactly two vertices have yet to be totally dominated, implying that the game will finish in at most two further moves. Thus, $\gtd(G) \le 3$. Consequently, $\gtd(G) = 3$. We show next that if $v$ is an arbitrary vertex of $G$, then $\gtd(G|v) = 2$. Let $v \in V(G)$ be arbitrary and consider the game played in $G|v$. By property~(b), there is a vertex $u$ of degree~$n-2$ in $G$ that is not adjacent to $v$. \D plays the vertex $u$ as his first move in $G|v$.  After Dominator's first move, all vertices of $G$, except for one, are totally dominated.  In particular, if $u\not= v$, then $u$ is not totally dominated.  If $u=v$, then the single vertex not in $N[v]$ is not totally dominated.  Thus, any (legal) move of \St in response to Dominator's first move finishes the game. Thus, $\gtd(G|v) \le 2$, implying that $G$ is a $3$-$\gtd$-critical graph.

Suppose secondly that $G$ is a $3$-$\gtd$-critical graph. By Corollary~\ref{c:co1}, property~(a) holds. By Lemma~\ref{l:lem2}, property~(b) holds. Hence it remains to show that property~(c) holds. For this purpose, let $v$ be an arbitrary vertex of $G$ of degree at most $n-3$, and consider the game played on $G|v$. By supposition, $\gtd(G|v) = 2$. Let $u$ be an optimal first move of \D in $G|v$. By Lemma~\ref{l:lem2B}, the vertex $u$ is not a neighbor of $v$. We show that $u$ has degree~$n-2$ in $G$. Suppose, to the contrary, that $u$ has degree at most~$n-3$. There are two cases to consider.

Suppose that $u \ne v$. In this case there exists a vertex $w$ of $G$ different from $u$ and $v$ that is not adjacent to~$u$. Every legal move of Staller in response to Dominator's first move $u$ finishes the game. In particular, we note that every neighbor of $u$ is a legal move for Staller, implying that every neighbor of $u$ is adjacent to~$w$. Thus, $N(u) \subseteq N(w)$. As observed earlier, $G$ is open twin-free, implying that $N(u) \subset N(w)$. Let $x \in N(w) \setminus N(u)$. If $x \ne v$, then \St plays the vertex $w$ on her first move, while if $x = v$, then \St plays the vertex $v$ on her first move. In both cases, the vertex $u$ is not yet totally dominated after the first two moves of the game, and so \St can force at least three moves to be made to finish the game, a contradiction.

Suppose that $u = v$. In this case there exists two distinct vertices, say $x$ and $y$, not in $N[u]$. Let $U = N(u)$ and let $\barU = V(G) \setminus N[u]$. We note that $\{x,y\} \subseteq \barU$. Let $S$ be the set of all legal moves of \St in response to Dominator's first move $u$. Since every legal move of \St finishes the game, we note that $\barU$ is an independent set in $G$ and every vertex in $S$ is adjacent to every vertex in $\barU$. In particular, this implies that every two distinct vertices in $\barU$ have the same open neighborhood, namely the set $S$. This contradicts the fact that $G$ is open twin-free. Hence, the vertex $u$ has degree~$n-2$ in $G$. Thus property~(c) holds.
\qed

Theorem~\ref{t:char3gtd} shows that there are infinitely many $3$-$\gtd$-critical graphs. The following result characterizes $3$-$\gtd$-critical graphs that can be obtained from the join of two graphs. Recall that the \emph{join} $G_1 + G_2$ of two graphs $G_1$ and $G_2$ is obtained from the disjoint union $G_1 \cup G_2$ of $G_1$ and $G_2$ by adding all possible edges joining $V(G_1)$ and $V(G_2)$. We say that a graph $G$ is the \emph{join} of two graphs if we can partition $V(G)$ into sets $V_1$ and $V_2$ so that $G = G_1 + G_2$, where $G_1 = G[V_1]$ and $G_2 = G[V_2]$.  The second main result of this section now reads as follows.

\begin{thm}
\label{t:join}
If $G_1$ and $G_2$ are vertex disjoint graphs, then $G_1 + G_2$ is $3$-$\gtd$-critical if and only if for each $i \in [2]$ one of the following holds. \\[-23.5pt]
\begin{enumerate}
\item The graph $G_i$ is $3$-$\gtd$-critical.
\item The graph $G_i \cong K_1 \cup K_{k}$ for some $k \ge 2$.
\end{enumerate}
\end{thm}
\proof Suppose that $G = G_1 + G_2$ is a $3$-$\gtd$-critical graph. If $u$ and $v$ are open twins in $G_1$, then $u$ and $v$ are open twins in $G$, contradicting Lemma~\ref{l:lem1}. Hence, neither $G_1$ nor $G_2$ has open twins. By Lemma~\ref{l:lem2}, the graph $G$ contains no dominating vertex, implying that neither $G_1$ nor $G_2$ contains a dominating vertex.

Suppose that $G_1$ has an isolated vertex $v$. Since $G$ has no open twins, the vertex $v$ is the only isolated vertex in $G_1$. Since $G_1$ has no dominating vertex, the graph $G_1$ therefore has at least three vertices. We show that in this case, $G_1 \cong K_1 \cup K_{k}$ for some $k \ge 2$.  Suppose, to the contrary, that there are two vertices $x$ and $y$ that are not adjacent in $G_1 - v$. By our earlier observations, both $x$ and $y$ have a neighbor in $G_1$. We now consider the graph $G|x$. Let $d_1$ be an optimal first move of \D in $G|x$. By Lemma~\ref{l:lem2B}, the vertex $d_1 \notin N_G(x)$, implying that $d_1$ is a vertex in $G_1$ that is not adjacent to $x$. Since $G$ is $3$-$\gtd$-critical, we note that $\gtd(G|x) = 2$, implying that every legal move for \St in $G|x$ in response to Dominator's first move $d_1$ completes the game.
If $d_1 = v$, then \St can play the vertex $y$ as her first move, and the game is not complete since the vertex $y$ is not yet totally dominated.
If $d_1 = x$, then \St can play a neighbor of $y$ in $G_1$ as her first move, and the game is not complete since the vertex $v$ is not yet totally dominated.
If $d_1 \ne v$ and $d_1 \ne x$,
then \St can play a neighbor of $d_1$ in $G_1$ as her first move, and the game is not complete since the vertex $v$ is not yet totally dominated. In all three cases, we have a contradiction. Therefore, $G_1 \cong K_1 \cup K_{k}$ for some $k \ge 2$, as claimed.

Suppose, next, that $G_1$ has no isolated vertex. In this case, we show that $G_1$ is $3$-$\gtd$-critical. As observed earlier, $G_1$ contains no dominating vertex and has no open twins. Consider the D-game played on $G_1$, let $v_1$ be an optimal first move of Dominator, and let $u_1$ be a vertex in $G_1$ not adjacent to $v_1$. Since $u_1$ and $v_1$ are not open twins, there is a vertex $w$ adjacent to exactly one $u_1$ and $v_1$. If \St plays the vertex $w$ on her first move, at least one vertex, namely $u_1$ or $v_1$, is not yet totally dominated after the first two moves of the game, implying that $\gtd(G_1) \ge 3$.

Let $v$ be an arbitrary vertex of $G_1$ and let $d_1$ be an optimal first move of \D in $G|v$. By Lemma~\ref{l:lem2B}, the vertex $d_1 \notin N_G(v)$, implying that $d_1$ is a vertex in $G_1$ that is not adjacent to $v$. Since $G$ is $3$-$\gtd$-critical, we note that $\gtd(G|v) = 2$, implying that every legal move for \St in $G|v$ in response to Dominator's first move $d_1$ completes the game.  In particular, every legal move for \St in $G_1|v$ in response to Dominator's first move $d_1$ completes the game played in $G_1|v$. Since $d_1$ has at least one neighbor in $G_1$ and every legal move for Staller in $G_1|v$ is a neighbor of $d_1$ in $G_1$, it follows that $\gtd(G_1|v) = 2$.

We have therefore shown in this case when $G_1$ has no isolated vertex that $\gtd(G_1) \ge 3$ and $\gtd(G_1|v) = 2$ for every vertex $v$ of $G_1$. Thus, the graph $G_1$ is $k$-$\gtd$-critical for some $k \ge 3$. It suffices for us to prove that $\gtd(G_1) = 3$. Let $v$ be a vertex of maximum degree in $G_1$, and let $d_1$ be an optimal first move of \D in $G_1|v$. As observed earlier, by Lemma~\ref{l:lem2B} the vertex $d_1 \notin N_G(v)$, implying that $d_1$ is a vertex in $G_1$ that is not adjacent to $v$.

We show that $d_1 = v$. Suppose, to the contrary, that $d_1 \ne v$. If $v$ has a neighbor that is not adjacent to $d_1$, then \St can play the vertex $v$ on her first move, resulting in at least one vertex, namely the vertex $d_1$, not yet totally dominated after the first two moves of the game, implying that $\gtd(G_1|v) \ge 3$, a contradiction. Therefore, every neighbor of $v$ is a neighbor of $d_1$ in $G_1$. However, by our choice of $v$ to be a vertex of maximum degree in $G_1$, this implies that $N(v) = N(d_1)$ in the graph $G_1$ (and the graph $G$), contradicting the fact that $G_1$ has no open twins. Hence, $d_1 = v$.

We now consider the game played in $G$, where \D plays as his first move the vertex $d_1$ ($=v$). Let $s_1$ be an optimal first move for \St in $G$ in response to Dominator's first move $d_1$. If $s_1$ is a legal move in $G_1|v$, then as observed earlier, $s_1$ is a neighbor of $d_1$ and her move completes the game in $G_1|v$, and therefore also completes the game played in $G_1$. If $s_1$ is not a legal move in $G_1|v$, then the only new vertex totally dominated by her move $s_1$ is the vertex $d_1$. Thus, after \St play her move $s_1$ all vertices in the set $N[v]$ are totally dominated. As observed earlier, every legal move in $G_1|N[v]$ completes the game. Thus, by playing the vertex $d_1$ as his first move in $G$, \D can guarantee that the game is completed in three moves or less.  Thus, $\gtd(G_1) \le 3$. As observed earlier, $\gtd(G_1) \ge 3$. Consequently, $\gtd(G_1) = 3$, implying that $G_1$ is a $3$-$\gtd$-critical graph.

Conversely, suppose that both graphs $G_1$ and $G_2$ satisfy (a) or (b) in the statement of the theorem. Let $G = G_1 + G_2$. If \D plays as his first move in the game played in $G$ a vertex from $G_1$ and as his second move a vertex from $G_2$, he can guarantee that the game will be complete in at most three moves, and so $\gtd(G) \le 3$. To show that $\gtd(G) \ge 3$, let $d_1$ be an optimal first move of \D played in $G$. Renaming $G_1$ and $G_2$, if necessary, we may assume that $d_1$ belongs to $G_1$. If $G_1$ is $3$-$\gtd$-critical, then \St can play as her first move in the game $G$ an optimal move for \St in the game played in $G_1$ with $d_1$ as Dominator's first move. If $G_1 \cong K_1 \cup K_{k}$ for some $k \ge 2$, then \St can play as her first move in the game $G$ a move in the clique $K_k$ different from $d_1$. In both cases, \St guarantees that at least three moves are needed to complete the game in $G$. Hence,  $\gtd(G) \ge 3$. Consequently, $\gtd(G) = 3$.

It remains for us to show that $\gtd(G|v) = 2$ for every vertex $v$ in $G$. Let $v$ be an arbitrary vertex of $G$. Renaming $G_1$ and $G_2$, if necessary, we may assume that $v$ belongs to~$G_1$. We consider the two possibilities for $G_1$ in turn.

Suppose firstly that $G_1$ is $3$-$\gtd$-critical and consider an optimal move $d_1$ for \D in the game played in $G_1|v$. \D now plays as his first move in the game played in $G|v$ the vertex $d_1$. Let $s_1$ be Staller's first move in $G|v$ in response to dominator's move $d_1$. If her move $s_1$ belongs to $G_2$, then the game is complete after her move since $\{d_1,s_1\}$ is a total dominating set of $G$. If her move $s_1$ belongs to $G_1$, then since $\gtd(G_1|v) = 2$ her move together with Dominator's first move $d_1$ totally dominate all vertices in $G_1$ except possibly for $v$, implying that $s_1$ and $d_1$ together totally dominate all vertices in $G$ except possibly for $v$. Thus, once again the game is complete after Staller's first move, implying that $\gtd(G|v) = 2$.

Suppose next that $G_1 \cong K_1 \cup K_{k}$ for some $k \ge 2$. If $v$ is the isolated vertex in $G_1$, then \D plays as his first move in the game played in $G|v$ an arbitrary vertex in the clique $K_k$ in $G_1$, thereby leaving exactly one vertex of $G$, namely the vertex played, not yet totally dominated in $G|v$ after his first move. Thus, every (legal) move of \St in response to Dominator's first move, completes the game. If $v$ belongs to the clique $K_k$ in $G_1$, then \D plays as his first move in the game played in $G|v$ the vertex $v$,  thereby leaving exactly one vertex of $G$ not yet totally dominated in $G|v$ after his first move, namely the vertex of $G$ that is isolated in $G_1$. Once again every (legal) move of \St in response to Dominator's first move, completes the game, implying that $\gtd(G|v) = 2$.
\qed

\section*{Acknowledgements}

Research of Michael Henning supported in part by the South African National Research Foundation and the University of Johannesburg. Sandi Klav\v zar acknowledges the financial support from the Slovenian Research Agency (research core funding No. P1-0297) and that the project (Combinatorial Problems with an Emphasis on Games, N1-0043) was financially supported by the Slovenian Research Agency. Supported by a grant from the Simons Foundation (Grant Number 209654 to Douglas F. Rall).

\end{document}